\documentclass{jnmp}

\usepackage{amsmath}

\JNMPnumberwithin{equation}{section}

\newtheorem{theorem}{Theorem}

\begin{document}

\renewcommand{\evenhead}{Z~Yin}
\renewcommand{\oddhead}{On the Cauchy Problem for a Nonlinearly Dispersive
Wave Equation}

\thispagestyle{empty}

\FirstPageHead{10}{1}{2003}{\pageref{yin-firstpage}--\pageref{yin-lastpage}}{Letter}

\copyrightnote{2003}{Z~Yin}

\Name{On the Cauchy Problem for a Nonlinearly Dispersive Wave
Equation} \label{yin-firstpage}

\Author{Zhaoyang YIN}

\Address{Permanent address: Department of Mathematics,\\
Zhongshan University, 510275 Guangzhou, China\\
E-mail: mcsyzy@zsu.edu.cn\\[10pt]
Present address: Department of Mathematics, Lund
University,\\
 P.O. Box 118, SE-221 00 Lund, Sweden\\
E-mail: yin@maths.lth.se}

\Date{Received May 15, 2002; Accepted August 19, 2002}

\begin{abstract}
\noindent We establish the local well-posedness for a new
nonlinearly dispersive wave equation and we show that the equation
has solutions that exist for indefinite times as well as solutions
which blowup in finite time. Furthermore, we derive an explosion
criterion for the equation and we give a sharp estimate from below
for the existence time of solutions with smooth initial data.
\end{abstract}

\section{Introduction}
In this letter we consider the Cauchy problem for the nonlinear
equation:
\begin{gather}
u_{t}-u_{txx}+2\omega u_{x}+3uu_{x} =
\gamma (2u_{x}u_{xx}+uu_{xxx}) ,\qquad t > 0,\quad x\in {\mathbb R}, \nonumber\\
u(0,x) = u_{0}(x),\qquad x\in {\mathbb R},
\end{gather}
with $\omega \geq 0$ and $\gamma \in {\mathbb R}$ fixed constants.

  With $\gamma =1$ in equation (1.1) we find the Camassa--Holm equation for
the unidirectional propagation of shallow water waves over a flat
bottom, $u(t,x)$ standing for the fluid velocity at time $t\geq 0$
in the spatial $x$ direction and $\omega$ being a nonnegative
parameter related to the critical shallow water speed
(see~\cite{C-H} and~\cite{J}).
 In this case equation~(1.1) has a
bi-Hamiltonian structure~\cite{F-F} and is completely integrable
(see~\cite{B-S-S,Cr,L}). The solitary waves are smooth if $\omega
>0$, they become peaked in the limiting case $\omega=0$, and
interact like solitons~\cite{C-H}.

For $\gamma = 0$, equation~(1.1) becomes the regularized long wave
equation,
  a well-known equation for surface
  waves in a channel~\cite{B-B-M}. All solutions are global and the solitary waves are
  smooth. Despite having a Hamiltonian
  structure, the equation is not integrable and its
  solitary waves are not solitons~\cite{D-E-G}.

For $\omega= 0$ and $\gamma \in {\mathbb R}$, equation~(1.1)
particularizes to a
  nonlinear dispersive model for finite-length and small-amplitude radial
  deformation waves in thin cylindrical compressible hyperelastic
  rods (see~\cite{Dai} and~\cite{D-H}), $u(t,x)$
representing the radial stretch relative to a~pre-stressed state.
The case $\gamma=1$ (limiting case of the Camassa--Holm equation)
has been extensively studied: the equation has global solutions
\cite{Cf, C-Ep} and also solutions which blow-up in finite time
\cite{Cf, C-E, C-Ep}. It is an integrable Hamiltonian system
\cite{B-S-S, Cr} and its solitary waves are peaked and interact
like solitons~\cite{C-H}.

   For $\omega > 0$ and $\gamma\in {\mathbb R}$, equation~(1.1)
seems not yet to have been discussed. The aim of this letter is to
prove local well-posedness of strong solutions to equation~(1.1)
for a~large class of initial data, to analyze global existence of
smooth solutions, as well as to get a blow-up criterion and a
sharp estimate from below  of the existence time for strong
solutions.

\section{Local Well-posedness}

To prove the local existence of solutions for equation~(1.1) we
will apply Kato's semigroup approach~\cite{K}. For this purpose,
let us reformulate equation~(1.1) as a quasi-linear evolution
equation for $y:=u-u_{xx}$,
\begin{equation}
 y_{t}+\gamma y_{x}u+2\gamma
yu_{x}+2\omega u_{x}+3(1-\gamma )uu_{x} = 0,\qquad t > 0,\quad
x\in {\mathbb R}.
\end{equation}
Note that if $p(x):=\frac{1}{2}e^{-|x|}$, $x\in {\mathbb R}$, then
$\left(1- \partial^{2}_{x}\right)f = p*f $ for all $f \in
L^{2}({\mathbb R})$ so $p \ast y=u $, where $\ast$ denotes
convolution. Using this identity, we can rewrite equation~(2.1) as
follows,
\begin{equation}
 u_{t}+\gamma uu_{x}+\partial_{x}p\ast
\left(\frac{3-\gamma}{2}u^{2}+\frac{\gamma}{2}u^{2}_{x}+2\omega
u\right) = 0,\qquad t > 0,\quad x\in {\mathbb R}.
\end{equation}
Equation~(2.2) is suitable for applying Kato's theory. A framework
that is similar to those provided in~\cite{C-Ep} and~\cite{Rb} can
be constructed. As an outcome, we obtain

\begin{theorem}\label{yin-theo1} Given $u_{0} \in H^{s}({\mathbb R})$,
$s>\frac{3}{2}$, there exists a maximal $T=T(\gamma, \omega,
u_{0})>0$, and a unique solution $u$ to equation~(1.1), such that
\[
u=u(\cdot,u_{0}) \in C\left([0,T);H^{s}({\mathbb R})\right)\cap
C^{1} \left([0,T);H^{s-1}({\mathbb R})\right).
\]
Moreover, the solution depends continuously on the initial data.
\end{theorem}

Theorem~\ref{yin-theo1} contains as particular cases (for various
choices of the constants $\gamma$ and $\omega$) the well-posedness
results obtained in \cite{C-Ep,H-M,L-O} and \cite{Rb}.

\section{Global solutions}

To answer affirmatively the question whether there are solutions
of equation~(1.1) that exist indefinitely in time, we prove that
for every choice of $\gamma \in {\mathbb R}$ and $\omega >0$ there
are smooth solitary wave solutions $u(t,x)=\varphi_c(x-ct)$
propagating at certain speeds $c>0$.

If $u(t,x)=\varphi_{c}(x-ct)$ is a smooth solitary wave solution
of equation~(1.1), then
\begin {equation}
-c\varphi'_{c}+c\varphi'''_{c}+3\varphi'_{c}\varphi_{c}
+2\omega\varphi'_{c}=2\gamma\varphi'_{c}\varphi''_{c}+
\gamma\varphi_{c}\varphi'''_{c}.
\end{equation}
Fix $x\in {\mathbb R}$ and integrate equation~(3.1) on
$(-\infty,x]$ to obtain
\begin {equation}
(2\omega-c)\varphi_{c}+c\varphi''_{c}+\frac{3}{2}(\varphi_{c})^{2}
=\frac{\gamma}{2}(\varphi'_{c})^{2}+\gamma\varphi_{c}\varphi''_{c}.
\end{equation}
Multiplying equation~(3.2) by $\varphi'_{c}$, an integration on
$(-\infty,x]$ leads to
\[
(\varphi'_{c})^{2}(c-\gamma
\varphi_{c})=\varphi_{c}^{2}(c-2\omega-\varphi_{c}).
\]
A detailed analysis of the above ordinary differential equation
yields the following result.

\begin{theorem} \label{yin-theo2}
Let $\omega > 0$ and $\gamma\in {\mathbb R}$. Then for every $c>0$
satisfying $c(\gamma -1)<2\omega\gamma$ and $c>2\omega$, there
exists a nontrivial smooth solitary wave solution
$u(t,x)=\varphi_{c}(x-ct)$ of equation~(1.1) and its profile
$\varphi_{c}$ decays exponentially to zero at infinity.
\end{theorem}

\section{Solutions with finite life-span}

In this section we address the question of the formation of
singularities for solutions to equation~(1.1). It turns out that
there are smooth initial data for which the corresponding solution
to equation~(1.1) does not exist globally in time. The precise
blowup scenario is described. We also give an explosion criterion
and a sharp estimate from below for the maximal existence time of
smooth solutions to equation~(1.1). Analogous to the case of the
Camassa--Holm equation~\cite{Ce}, equation~(1.1) can be written in
Hamiltonian form and has the invariant
 \begin{equation}
   E(u)=\int_{{\mathbb R}}\left(u^{2}+u^{2}_{x}\right)dx .
\end{equation}
The invariance of $E(u)$ ensures that all solutions of
equation~(1.1) are uniformly bounded as long as they exist.
Arguments similar to those developed in~\cite{C-S2} for the
elasticity model can be pursued to obtain the following result.

\begin{theorem}\label{yin-theo3}
Given $u_{0} \in H^{s}({\mathbb R})$, $s > \frac{3}{2}$, the
solution $u=u(\cdot,u_{0}$) of equation~(1.1) is uniformly
bounded. Blow up in finite time $T < + \infty$ occurs if and only
if
\[
\liminf _{t \uparrow T}\, \{\inf _{x \in {\mathbb R}}[\,\gamma
u_{x} (t,x)] \} = - \infty .\]
\end{theorem}

 As an interesting outcome of
the previous result, note that in the special case $\gamma =0$,
all solutions to equation~(1.1) are global. Also, in the
particular cases when equation~(1.1) models water waves or waves
on an elastic rod, the previously described blowup corresponds
physically to a breaking wave, respectively rod.

We discuss now the question of finite time blow-up of solution to
equation~(1.1) with rather general initial data. Let $T>0$ be the
maximal existence time of the solution $u(t,\cdot)$ of
equation~(1.1) with initial data $u_0 \in H^s({\mathbb R})$, $s >
\frac{3}{2}$. Differentiate equation~(2.2) with respect to $x$,
and multiply the obtained equation by $\gamma$. Since
$\partial_{x}^{2}p\ast f=p\ast f-f$, we have
\begin{gather}
\gamma u_{tx}=-\frac{\gamma^{2}}{2}u_{x}^{2}-\gamma^{2}uu_{xx}
+\frac{(3-\gamma)\gamma}{2}u^{2}\nonumber\\
\phantom{\gamma u_{tx}=}{}+2\omega\gamma u-p\ast
\left(\frac{(3-\gamma)\gamma}{2}u^{2}+\frac{\gamma^{2}}{2}u^{2}_{x}+2\omega\gamma
u\right).
\end{gather}
Define $m(t)=\gamma u_{x}(t,\xi(t))=\inf\limits_{x\in {\mathbb
R}}\{\gamma u_{x}(t,x)\}$, where $\xi(t)$ is any point where this
minimum is attained. Since we deal with a minimum,
$u_{xx}(t,\xi(t))=0$ for all $t\in [0,T)$. On the other hand, the
regularity of the solution $u(t,x)$ provided by
Theorem~\ref{yin-theo1} allows us to infer from the main result
in~\cite{C-E} that the function $m(t)$ is almost everywhere
differentiable on $(0,T)$, with
\[
\frac{dm}{dt}=\gamma u_{tx}(t,\xi(t))    \qquad \mbox{a.e.} \quad
\mbox{on}\ (0,T).
\]
Hence we obtain that
\begin{gather}
m'(t)=-\frac{1}{2}m^{2}(t)+\frac{(3-\gamma)\gamma}{2}u^{2}(t,\xi(t))+2\omega\gamma
u(t,\xi(t))\nonumber\\
\phantom{m'(t)=}{}-p\ast
\left(\frac{(3-\gamma)\gamma}{2}u^{2}+\frac{\gamma^{2}}{2}u^{2}_{x}+2\omega\gamma
u\right)(t,\xi(t)).
\end{gather}
Involved manipulations of the equations (4.2) and (4.3) and
estimates analogous to those in \cite{D} for the special case of
the Camassa--Holm equation lead us to the following result.

\begin{theorem}\label{yin-theo4}  (i) Given $u_{0} \in H^{s}({\mathbb R})$, $s >
\frac{3}{2}$, and assume that we can find $x_{0} \in {\mathbb R}$
with
\[
\gamma u'_{0}(x_{0}) <
 - \left(\frac{|(\gamma -3)\gamma|}{2}E(u_{0})+
4\sqrt{2} \omega|\gamma| (E(u_{0}))^{1/2}\right) ^{1/2}, \qquad
\gamma \neq 0 .
\]
Then the corresponding solution to equation~(1.1) blows up in
finite time.

(ii) Given $u_{0} \in H^{s}({\mathbb R})$, $s > \frac{3}{2}$, the
maximal existence time of the solution to equation~(1.1) with
initial data $u_0$ is estimated from below by
\[
T_{u_{0}}=\left\{\begin{array}{ll} \displaystyle
\frac{-2\arctan\frac{\left(\frac{(3-\gamma)\gamma}{2}E(u_{0})
+4\sqrt{2}\omega|\gamma|(E(u_{0}))^{1/2}\right)^{1/2}}{m(0)}}
{\left(\frac{(3-\gamma)\gamma}{2}E(u_{0})+4\sqrt{2}\omega|\gamma|
(E(u_{0}))^{1/2}\right)^{1/2}},&\displaystyle if\quad
0<\gamma<\frac{3}{2},\vspace{1mm}\\
\displaystyle
\frac{-2\arctan\frac{\left(\frac{\gamma^{2}}{2}E(u_{0})+4\sqrt{2}\omega|\gamma
|(E(u_{0}))^{\frac{1}{2}}\right)^{1/2}}{m(0)}}{\left(\frac{\gamma^{2}}{2}E(u_{0})
+4\sqrt{2}\omega|\gamma|(E(u_{0}))^{1/2}\right)^{1/2}},&\displaystyle
if\quad
\frac{3}{2}\leq\gamma\leq3,\vspace{1mm}\\
\displaystyle
\frac{-2\arctan\frac{\left(\frac{(2\gamma-3)\gamma}{2}E(u_{0})
+4\sqrt{2}\omega|\gamma
|(E(u_{0}))^{1/2}\right)^{1/2}}{m(0)}}{\left(\frac{(2\gamma-3)\gamma}{2}E(u_{0})
+4\sqrt{2}\omega|\gamma
|(E(u_{0}))^{1/2}\right)^{1/2}},&\displaystyle if\quad \gamma> 3\
\  or\ \ \gamma < 0,
\end{array}\right.
\]
where $ m(0):=\inf\limits_{x\in {\mathbb R}}\,\{{\gamma
\partial_xu_{0} (x)} \}$. The previous estimate is sharp.
\end{theorem}

\subsection*{Acknowledgments}

This work was performed as a Visiting Researcher
 at Lund University. The author is very pleased to acknowledge the support
 and encouragement of Professor A~Constantin as this work has developed.
 The author is also grateful to Professor J~L~Bona for his interest in
 these investigations and for helpful discussions. This work was partially
 supported by the National Natural Science Foundation of China and the
 Foundation of Zhongshan University Advanced Research Center.

\label{yin-lastpage}

\end{document}